\documentclass[11pt]{amsart}
\usepackage{amsmath, amsfonts, amscd, epsfig, amssymb}

\newcommand{\RR}{{\mathbb R}}

\newcommand\adots{\mathinner{\mkern2mu\raise1pt\hbox{.}
\mkern3mu\raise4pt\hbox{.}\mkern1mu\raise7pt\hbox{.}}}

\renewcommand{\div}{{\rm div}}
\newcommand{\curl}{{\rm curl}}

\newtheorem{theo}{Theorem}[section]

\newtheorem{cor}[theo]{Corollary}
\newtheorem{lem}[theo]{Lemma}
\newtheorem{defi}[theo]{Definition}
\newtheorem{ass}[theo]{Assumption}

\newtheorem{rem}[theo]{Remark}

\numberwithin{equation}{section}

\begin{document}

\title{Structural conditions for full MHD equations}

\author{Bongsuk Kwon}
\address{Indiana University, Bloomington, IN 47405, USA}
\email{bkwon@indiana.edu}
\thanks{Thanks to Kevin Zumbrun for suggesting the problem. This work was
partially supported under NSF grant DMS-0300487}

\begin{abstract}
In this paper, we investigate the characteristic structure of the
full equations of magnetohydrodynamics (MHD) and show that it
satisfies the hypotheses of a general variable-multiplicity
stability framework introduced by M\'{e}tivier and Zumbrun, thereby
extending to the general case various results obtained by
M\'{e}tivier and Zumbrun for the isentropic equations of MHD.

\end{abstract}

\maketitle \thispagestyle{empty}



\section{Introduction}

Stability of hyperbolic boundary value problems and the related
shock problems has been studied by a number of authors for years.
There are satisfactory theories for symmetric systems with
dissipative boundary conditions or hyperbolic systems with constant
multiplicities provided that the boundary conditions satisfy
Lopatinki conditions. Many applications fall into one of these two
cases. But the equations of magnetohydrodynamics (MHD) fail to
satisfy the assumption of constant multiplicities. Moreover, though
they are symmetrizable, boundary conditions arising in the shock
stability problem are not dissipative. This motivates the further
study of the variable multiplicities case. It is proved by
\cite{MeZ} for symmetrizable systems that under certain conditions
on eigenvalues, the uniform Lopatinski condition is equivalent to
the maximal energy estimate which gives linearized and nonlinear
stability. Specifically, they extend the construction of Kreiss'
symmetrizers to the case that eigenvalues are not geometrically
regular provided that they are ``totally nonglancing" or ``linearly
splitting" in a sense that they define. (See section 3 for their
definitions) This yields the maximal estimate when the uniform
Lopatinski condition is satisfied. In \cite{MeZ}, it is also showed
that all characteristics of isentropic MHD are either
``geometrically regular" or ``totally nonglancing". Hence, a
satisfactory framework for the study of linearized and nonlinear
stability of isentropic MHD is achieved, through which various
results are obtained. The goal of this paper is to investigate the
eigenvalues of full MHD equations. It turns out that all
characteristics are algebraically regular. More specifically, they
are either geometrically regular or totally nonglancing as the
isentropic case. As a consequence, all results obtained in
\cite{MeZ} for the isentropic case extend immediately to the full
MHD equations.

\section{Definitions}
Consider a first-order system
\begin{equation}\label{1}
L(p,\partial_t,\partial_{x})U=U_t+\sum_{j=1}^{d} A^{j}(p)U_{x_j},
\end{equation}
and its symbolic operator
\begin{equation}\label{2}
\widetilde{L}(p,\tau,\xi)=\tau Id + \sum_{j=1}^{d}\xi_j A^{j} (p).\\
\end{equation}
The characteristic polynomial of \eqref{2} is $P(p,\tau,\xi):= \det 
\widetilde{L}(p,\tau,\xi)$.
%

\begin{defi}
\smallbreak$(\mathcal{D}1)$ The system \eqref{1} is hyperbolic if $\sum_{j=1}^{d}\xi_j A_j (p)$ has all real, semisimple
eigenvalues for all $\xi\in\RR^{d}$.
\smallbreak$(\mathcal{D}2)$ The system \eqref{1} is symmetric hyperbolic in the
sense of Friedrichs if there is a positive definite smooth
symmetrizer $S(p)$ such that $S(p)A_j(p)$ is symmetric for all
$j=1,...,d$.
\end{defi}

Let $(p_0,\tau_0,\xi_0)$ be a root of $P(p,\tau,\xi)=0$ whose
multiplicity is $m$ in $\tau$.
\begin{defi}
 $(a)$ $(p_0,\tau_0,\xi_0)$ is algebraically regular if on a
neighborhood
 $\omega$ of $(p_0,\xi_0)$, there are $m$ smooth real functions
 $\lambda_j(p,\xi)$, analytic in $\xi$, such that for $(p,\xi)$,
 \begin{equation}\label{3}
 P(p,\tau,\xi)=e(p,\tau,\xi) \Pi_{j=1}^{m} (\tau+\lambda_j(p,\xi))
\end{equation}
where $e$ is a polynomial in $\tau$ and smooth coefficients in $\xi$
and $p$ such that $e(p_0,\tau_0,\xi_0)\neq0$.
\smallbreak $(b)$ $(p_0,\tau_0,\xi_0)$ is geometrically regular if in addition
 there are $m$ smooth vectors $v_j(p,\xi)$, analytic in $\xi$, such
 that
 \begin{equation}\label{5}
 A(p,\xi)v_j(p,\xi)=\lambda_j(p,\xi)v_j(p,\xi)
 \end{equation}
 where $v_1,...v_m$ are linearly independent.
\end{defi}

\begin{defi}
A root $(p_0,\tau_0,\eta_0,\xi_0)$ of $P$, of multiplicity $m$ in
$\tau$ is nonglancing if the $m$ th order of Taylor expansion of $P$
at $(p_0,\tau_0,\eta_0,\xi_0)$ satisfies
\begin{equation} \label{taylor}
P^{m}(dx)\neq0
\end{equation}
where $dx$ is the (space-time) conormal to the boundary. Moreover,
it is called totally nonglancing if it is either totally incoming or
outgoing.
\end{defi}

\begin{ass}\label{ass1} Suppose that $L$ is symmetric hyperbolic and
noncharateristic with respect to the boundary $x_d=0$, and that its
characteristic roots are geometrically regular or totally
nonglancing.
\end{ass}

Consider a hyperbolic system
\begin{equation}\label{5}
L(p,\partial_t,\partial_{\tilde{x}},\partial_{x_d})=A_d(p)(\partial_{x_d}+G(p,\partial_t,\partial_{\tilde{x}}))
\end{equation}
with the assumption that the boundary $\{x_d=0\}$ is
noncharacteristic, i.e., \text{ det }$A_d(p)\neq 0$.

The classical plane wave analysis yields the boundary value problems
\begin{equation}\label{4}
\begin{aligned}
\left\{%
\begin{array}{l}
 L(p,\zeta) \tilde{U}= \partial_{x_d} \tilde{U}+ i G(p, \zeta) \tilde{U}=f,\\
 M(p,\zeta) \tilde{U} _ {|_{x_d =0}}=g,\\
\end{array}%
\right.
\end{aligned}
\end{equation}
where \begin{equation} G(p,\zeta)=A_d^{-1}(p)((\tau-i\gamma) Id
+\sum_{j=1}^{d-1} \eta_j A^{j}(p))\\
\end{equation}
with a new variable $\zeta=(\tau-i\gamma,\eta)$.\\

Let $\textbf{E}_{-}(p,\zeta)$ be the invariant subspace of
$G(p,\zeta)$ associated to the eigenvalues in $\{\text{Im } \mu <
0\}$.

\begin{defi}
The Lopatinski determinant associated with $(L,M)(p,\zeta)$ is defined by
\begin{equation}\label{6}
D(p,\zeta):= \det (\textbf{E}_{-}(p,\zeta), \ker M(p,\zeta)).
\end{equation}

We say that $(L,M)(p,\zeta)$ satisfies the uniform Lopatinski condition on a
neighborhood $\omega$ if there is a constant $c>0$ such that

\begin{equation}\label{7}
|D(p,\zeta)|\geqslant c \text{ for } \forall (p,\zeta) \in \omega
\times S_+^{d}.
\end{equation}
\end{defi}


\section{Computation}

The inviscid full MHD Equations are
\begin{equation}\label{firstorder}
\begin{aligned}
\left\{%
\begin{array}{l}
 \vspace{2mm} \displaystyle \rho_t + \div(\rho u) = 0,\\
\vspace{2mm} \displaystyle  \rho(u_t + u\cdot\nabla u) + \nabla P - \curl B \times B=0,\\
\vspace{2mm} \displaystyle  \rho e_\theta(\theta_t+u\cdot\nabla\theta) +\theta P_\theta \div u= 0, \\
 \displaystyle  B_t+\curl(B\times u)=0
\end{array}%
\right.
\end{aligned}
\end{equation}
with the constraint
\begin{equation}\label{secondorder}
\div B = 0,
\end{equation}
where $\rho\in  \RR$ is density, $u\in \RR^3$ is
velocity, $P(\rho,\theta)\in  \RR $ is pressure,
$B\in\RR^{3}$ is magnetic field and $\theta \in\RR$ is
temperature.\\
Using the identity $\curl(B\times u)=
(\text{div}u)B+u\cdot\nabla B -
B\cdot\nabla u -(\div B)u $ and the constraint \eqref{secondorder}, \eqref{firstorder} can be put in the symmetric form

\begin{equation}\label{thirdorder}
\begin{aligned}
\left\{%
\begin{array}{l}
\vspace{2mm} \displaystyle  \rho_t + u\cdot \nabla \rho + \rho \div u =0,\\
\vspace{2mm} \displaystyle  u_t+u\cdot\nabla u + \rho^{-1} (P_\rho\nabla\rho+P_\theta\nabla\theta ) - \rho^{-1}(\curl B\times B)=0,\\
 \vspace{2mm} \displaystyle \theta_t + u\cdot \nabla\theta + \frac{\theta P_\theta}{\rho
  e_\theta} \div u=0,\\
 \displaystyle B_t + u\cdot\nabla B + (\div u) B -B\cdot\nabla u=0.
\end{array}%
\right.
\end{aligned}
\end{equation}
Linearizing \eqref{thirdorder} about $(\rho,u,\theta,B)$, we have
the linearized equations
\begin{equation}\label{forthorder}
\begin{aligned}
\left\{%
\begin{array}{l}
\vspace{2mm} \displaystyle  D_t \dot{\rho} + \rho \div\dot{u}=0,\\
 \vspace{2mm} \displaystyle D_t \dot{u}+\rho^{-1}(P_\rho\nabla\dot{\rho}+P_\theta\nabla\dot{\theta})-\rho^{-1}(\curl\dot{B}\times
  B)=0,\\
 \vspace{2mm} \displaystyle D_t \dot{\theta}+\frac{\theta P_\theta}{\rho e_\theta} \div \dot{u}=0,\\
 \displaystyle D_t \dot{B}+(\div \dot{u})B-B\cdot\nabla\dot{u}=0
\end{array}%
\right.
\end{aligned}
\end{equation}
where $D_t=\partial_t+u\cdot\nabla$ is the material derivative.

By taking Fourier Transform in $x$ variable, we have the symbol
matrix associated to \eqref{forthorder}
\begin{align}
\tilde{A}(U,\xi)&:=(-u\cdot\xi)I+\sum_{j} A^j(U)\xi_j \\
&=\left(%
\begin{array}{cccc}
  0 & \rho \xi & 0 & 0 \\
  \rho^{-1} P_\rho \xi^{tr} & 0 & \rho^{-1} P_\theta \xi^{tr} & \rho^{-1}(\xi^{tr}B-(B\cdot\xi)I)\\
  0 & \rho^{-1} P_\theta \theta/e_\theta \xi & 0 &0\\
  0 & B^{tr}\xi-(B \cdot \xi)I & 0 & 0\\
\end{array}\nonumber%
\right).
\end{align}

So, we can express the linear part of \eqref{forthorder} in the
matrix form
\begin{align}
\tilde{A}(U,\xi) \left(%
\begin{array}{cccc}
\dot{\rho}\\
\dot{u}\\
\dot{\theta}\\
\dot{B}\\
\end{array}%
\right)=
\left(%
\begin{array}{cccc}
\rho(\xi\cdot \dot{u})\\
\rho^{-1}P_\rho\dot{\rho}\xi^{tr}+\rho^{-1}P_\theta\dot{\theta}\xi^{tr}+\rho^{-1}(B \cdot \dot{B}) \xi^{t}-\rho^{-1}(B\cdot \xi)\dot{B}^{t}\\
\rho^{-1}P_\theta \theta /e_\theta(\xi\cdot\dot{u})\\
B^{t}(\xi\cdot\dot{u})-(B\cdot\xi)\dot{u}^{t}\\
\end{array}%
\right).
\end{align}

Now, consider the eigenvalue equation associated to
\eqref{forthorder}
\begin{align}
\tilde{A}(U,\xi) \left(%
\begin{array}{cccc}
\dot{\rho}\\
\dot{u}\\
\dot{\theta}\\
\dot{B}\\
\end{array}%
\right)=
\lambda \left(%
\begin{array}{cccc}
\dot{\rho}\\
\dot{u}\\
\dot{\theta}\\
\dot{B}\\
\end{array}%
\right).
\end{align}

Equivalently, we have

\begin{align}\label{eveq}
\left\{%
\begin{array}{l}
\vspace{2mm} \displaystyle  \rho(\xi\cdot\dot{u})=\lambda \dot{\rho},\\
\vspace{2mm} \displaystyle  \rho^{-1}P_\rho\dot{\rho} \xi^{t}+\rho^{-1}P_\theta\dot{\theta} \xi^{t}+\rho^{-1}(B \cdot \dot{B}) \xi^{t}-\rho^{-1}(B\cdot \xi)\dot{B}^{t}=\lambda\dot{u}^{t},\\
\vspace{2mm} \displaystyle  \rho^{-1}P_\theta \theta /e_\theta(\xi\cdot\dot{u})=\lambda\dot{\theta},\\
\displaystyle (\xi\cdot\dot{u})B^{t}-(B\cdot\xi)\dot{u}^{t}=\lambda\dot{B}.\\
\end{array}%
\right.
\end{align}

We adopt the variable $\xi$ for the spatial frequencies and let
\begin{align}
\xi &=|\xi|\hat{\xi}, \hspace{3mm}u_\|=\hat{\xi}\cdot u, \hspace{3mm} u_\bot=u-u_\|\hat{\xi}.
\end{align}

Then we can rewrite the eigenvalue equations \eqref{eveq} as
\begin{align}\label{char}
\left\{%
\begin{array}{l}
\vspace{2mm} \displaystyle \rho\dot{u}_\|=\tilde{\lambda}\dot{\rho},\\
\vspace{2mm} \displaystyle  P_\rho\dot{\rho}+P_\theta\dot{\theta}+ B_\bot \cdot \dot{B}_\bot=\tilde{\lambda}\rho\dot{u}_\|,\\
\vspace{2mm} \displaystyle  -\mu_0^{-1}B_\|\dot{B}_\bot=\tilde{\lambda}\rho\dot{u}_\bot,\\
\vspace{2mm} \displaystyle  P_\theta\theta/e_\theta\dot{u}_\|=\rho\tilde{\lambda}\dot{\theta},\\
 \vspace{2mm} \displaystyle 0=\tilde{\lambda}\dot{B}_\|,\\
\displaystyle  \dot{u}_\|B_\bot-B_\|\dot{u}_\bot=\tilde{\lambda}\dot{B}_\bot.
\end{array}%
\right.
\end{align}

 Let $\dot{\sigma}:=\dot{\rho}/\rho$ and
$\dot{v}:=\dot{B}/\sqrt{\rho}$. The characteristic system
 \eqref{char} reads
\begin{align}
\left\{%
\begin{array}{l}
\vspace{2mm} \displaystyle  \dot{u}_\|=\tilde{\lambda}\dot{\sigma},\\
\vspace{2mm} \displaystyle  P_\rho\dot{\sigma}+\rho^{-1}P_\theta\dot{\theta}+v_\bot \cdot \dot{v}_\bot=\tilde{\lambda}\dot{u}_\|,\\
 \vspace{2mm} \displaystyle -v_\|\dot{v}_\bot=\tilde{\lambda}\dot{u}_\bot,\\
\vspace{2mm} \displaystyle  \frac{P_\theta\theta}{e_\theta\rho}\dot{u}_\|=\tilde{\lambda}\dot{\theta},\\
 \vspace{2mm} \displaystyle 0=\tilde{\lambda}\dot{v}_\|,\\
 \displaystyle v_\bot\dot{u}_\|-v_\|\dot{u}_\bot=\tilde{\lambda}\dot{v}_\bot.
\end{array}%
\right.
\end{align}

Now, we introduce new variables
\begin{equation}\label{key}
\begin{aligned}
\left\{%
\begin{array}{l}
\vspace{2mm} \displaystyle \dot{x}:=P_\rho\dot{\sigma}+\rho^{-1}P_\theta\dot{\theta},\\
\dot{y}:=P_\theta\theta \dot{\sigma}-e_\theta\rho\dot{\theta}.\\
\end{array}%
\right.
\end{aligned}
\end{equation}

 Note that the transformation
$(\dot{\sigma},\dot{\theta})\mapsto(\dot{x},\dot{y})$ is invertible
if and only if $\rho P_\rho e_\theta+\rho^{-1}\theta
P_\theta^{2}\neq0$. So, it is invertible since
$P_\rho+\frac{P_\theta^{2}\theta}{\rho^{2}e_\theta}>0$. Using new
variables $(\dot{x}, \dot{y})$, we have the equivalent system
\begin{equation}\label{order}
\begin{aligned}
\left\{%
\begin{array}{l}
 \vspace{2mm} \displaystyle (P_\rho+\frac{P_\theta^{2}\theta}{\rho^{2}e_\theta})\dot{u}_\|=\tilde{\lambda}\dot{x},\\
\vspace{2mm} \displaystyle  0=\tilde{\lambda}\dot{y},\\
\vspace{2mm} \displaystyle  \dot{x}+v_\bot \cdot \dot{v}_\bot=\tilde{\lambda}\dot{u}_\|,\\
\vspace{2mm} \displaystyle  -v_\|\dot{v}_\bot=\tilde{\lambda}\dot{u}_\bot,\\
 \vspace{2mm} \displaystyle 0=\tilde{\lambda}\dot{v}_\|,\\
  v_\bot\dot{u}_\|-v_\|\dot{u}_\bot=\tilde{\lambda}\dot{v}_\bot.\\
\end{array}%
\right.
\end{aligned}
\end{equation}

Note that the zero-eigenspace of $\tilde{A}$ is $ \mathbb{E}_0
:=\{\dot{x}=0, \dot{u}=0, \dot{v}_\bot=0\} $, which is of dimension
2. On the space $\mathbb{E}_0$, $\tilde{A}=0$. This implies
$A=u\cdot\xi$ on $\mathbb{E}_0$. Thus, $u\cdot\xi$ is an eigenvalue
of $A$ of constant multiplicity 2. By inspection, it is both
algebraically and geometrically regular.

Taking a basis of $\xi^{\bot}$ such that $v_\bot=(b,0)$ and letting
$a=v_\|$, we can write \eqref{order} in the coordinates of
$(\dot{x},\dot{u}_\|,
\dot{u}_\bot^{1},\dot{u}_\bot^{2},\dot{v}_\bot^{1},\dot{v}_\bot^{2},\dot{y},\dot{v}_\|)^{t}$
as follows:
\begin{equation}\label{11}
-\tilde{A_0}+\tilde{\lambda}I=
\begin{pmatrix}
\tilde{\lambda}&-\gamma^{2}&0&0&0&0&0&0\\
-1&\tilde{\lambda}&0&0&-b&0&0&0\\
0&0&\tilde{\lambda}&0&a&0&0&0\\
0&0&0&\tilde{\lambda}&0&a&0&0\\
0&-b&a&0&\tilde{\lambda}&0&0&0\\
0&0&0&a&0&\tilde{\lambda}&0&0\\
0&0&0&0&0&0&\tilde{\lambda}&0\\
0&0&0&0&0&0&0&\tilde{\lambda}\\
\end{pmatrix}=0
\end{equation}
where
$\gamma^{2}=P_\rho+\frac{P_\theta^{2}\theta}{\rho^{2}e_\theta}>0$.\\

So, the associated characteristic polynomial $P(x)$ is
\begin{align}
P(x)&=x^{2}(x^{2}-a^{2})\big((x^{2}-a^{2})(x^{2}-(P_\rho+\frac{P_\theta^{2}\theta}{\rho^{2}e_\theta}))-b^{2}x^{2}\big)\\
&=x^{2}(x^{2}-a^{2})\big((x^{2}-a^{2})(x^{2}-\gamma^{2})-b^{2}x^{2}\big).\nonumber
\end{align}
Thus,

\begin{equation}\label{evalues}
\begin{aligned}
\left\{%
\begin{array}{l}
\vspace{2mm} \displaystyle  \tilde{\lambda}^{2}=0,\\
\vspace{2mm} \displaystyle  \tilde{\lambda}^{2}=a^{2},\\
 \vspace{2mm} \displaystyle \tilde{\lambda}^{2}=c_f^{2}:=\frac{1}{2}\big((\gamma^{2}+h^{2})+\sqrt{(\gamma^{2}-h^{2})^{2}+4b^{2}\gamma^{2}}\big),\\
\displaystyle  \tilde{\lambda}^{2}=c_s^{2}:=\frac{1}{2}\big((\gamma^{2}+h^{2})-\sqrt{(\gamma^{2}-h^{2})^{2}+4b^{2}\gamma^{2}}\big)\\
\end{array}%
\right.
\end{aligned}
\end{equation}
where $h^{2}=a^{2}+b^{2}=|B|^{2}/\rho$.

\begin{rem}
Coordinate change \eqref{key} is the key step in the computation,
leading to block decomposition \eqref{11}.
\end{rem}
\section{Results}
In the analysis of eigenvalues of \eqref{firstorder}, what is left
is the same as the corresponding steps carried out in \cite{MeZ} for
the isentropic case.
Specifically, in the isentropic case treated in \cite{MeZ}, the
matrix corresponding to \eqref{11} has the same form, with identical
$6\times 6$ upper lefthand block and lower lefthand block of the
same form $\tilde \lambda I_k$, but with dimension $k=1$ instead of
$k=2$ as in the present case. Thus, decomposing into these two
blocks, all the analysis goes as before.
We summarize for completeness.

\begin{lem} Let
$\gamma^{2}=P_\rho+\frac{P_\theta^{2}\theta}{\rho^{2}e_\theta}>0$.
The 8 eigenvalues are
\begin{align}\label{lemmma1}
\left\{%
\begin{array}{l}
\vspace{2mm} \displaystyle  \lambda_{0}=\xi\cdot u  \text{ (of multiplicity $2$)},\\
 \vspace{2mm} \displaystyle \lambda_{\pm 1}=\xi\cdot u\pm c_{s}|\xi|, \\
\vspace{2mm} \displaystyle  \lambda_{\pm 2}=\xi\cdot u\pm (\xi\cdot B)/\sqrt{\rho},\\
  \lambda_{\pm 3}=\xi\cdot u\pm c_{f}|\xi|\\
\end{array}%
\right.
\end{align}
where
\begin{equation}
c_{f}=\frac{1}{2}((\gamma^{2}+h^{2})+\sqrt{(\gamma^{2}-h^{2})^{2}+4b^{2}\gamma^{2}}),
\end{equation}
\begin{equation}
c_{s}=\frac{1}{2}((\gamma^{2}+h^{2})-\sqrt{(\gamma^{2}-h^{2})^{2}+4b^{2}\gamma^{2}}),
\end{equation}
\begin{equation} h^{2}=|B|^{2}/\rho,
\end{equation} and
\begin{equation}
b^{2}=|\hat{\xi}\times B|^{2}/\rho.
\end{equation}
\end{lem}

\begin{lem} Assume that $0<|B|^{2}\neq \rho\gamma^{2}$ where $
\gamma^{2}=P_\rho+\frac{P_\theta^{2}\theta}{\rho^{2}e_\theta}>0$. \\
(a) When $\xi\cdot B\neq 0$ and $\xi\times B\neq 0$, the eigenvalues
$\lambda_{\pm 1}, \lambda_{\pm 2}$ and $\lambda_{\pm 3}$ are simple
and $ \lambda_{0}$ is semisimple of constant multiplicity 2 (geometrically regular). \\
(b) On the manifold $\xi \cdot B=0, \xi\neq 0$, the eigenvalues
$\lambda_{\pm 3}$ are simple and the multiple eigenvalues
$\lambda_{0}=\lambda_{\pm 1}=\lambda_{\pm 2}$ are geometrically
regular.\\
(c) On the manifold $\xi\times B=0, \xi\neq 0$, the multiple
eigenvalue $\lambda_{0}$ is semisimple. When $|B|^{2}<\rho\gamma^{2}
$ (resp. $ |B|^{2}>\rho\gamma^{2}), \lambda_{\pm 3} (\text{resp.
}\lambda_{\pm 1}$) are simple and
$\lambda_{+2}=\lambda_{+1}(\text{resp.
}\lambda_{+3})\neq\lambda_{-2}=\lambda_{-1}(\text{resp.
}\lambda_{-3})$ are double and algebraically regular, and not
geometrically regular (totally nonglancing, provided that
$u_3-\sigma\neq\pm B_3/\sqrt{\rho}$).
\end{lem}
\begin{proof}
The characteristic polynomial is the same as the one for the
isentropic case, except additional factor $x$. The zero-eigenspace
is of dimension 2, independent of $\xi$. Thus, the eigenvalue
$\lambda_0=\xi\cdot u$ is geometrically regular. For other
eigenvalues, the result follows from the calculation in \cite{MeZ}.
\end{proof}

\begin{cor} For the linearized equations of (full) MHD, the uniform Lopatinski
condition is equivalent to the uniform stability estimate (implying
linearized and nonlinear stability).
\end{cor}

 \begin{cor}
\label{mhdlop} Under the generically satisfied conditions $0 <\vert
B^\pm \vert^2 \ne \rho^\pm (\gamma^\pm)^2$,
 the uniform Lopatinski condition implies
linearized and nonlinear stability of noncharacteristic Lax-type MHD
shocks.
\end{cor}

As an application, we recover for the full equations of MHD the
result of Blokhin and Trakhinin that MHD shocks are stable in the
small-magnetic field limit $B\to 0$.

\begin{cor}\label{pert}
In the $B\to 0$ limit, Lax-type MHD shocks approaching a
noncharacteristic limiting fluid-dynamical shock 
satisfy the Lopatinski condition if and only if it is
satisfied by the limiting fluid-dynamical shock (i.e., the shock
satisfies the physical stability conditions of Erpenbeck--Majda
\cite{Er, Maj}, in which case the MHD shocks are linearly and
nonlinearly stable. In particular, for an ideal gas equation of
state, they are always stable in the $B\to 0$ limit.
\end{cor}

\begin{rem}
\textup{ Verification of the uniform Lopatinski condition has been
carried out numerically or analytically for several interesting
cases; see \cite{BT.1, BT.2}.
}
\end{rem}

\end{document}